\documentclass[11pt]{llncs}
\usepackage{epsfig}
\usepackage{graphicx}
\usepackage{color}
\usepackage{amsmath,amssymb}    
\usepackage{latexsym}
\usepackage{hhline}
\spnewtheorem{observation}{Observation}{\bfseries}{\itshape}
\spnewtheorem{comment}{Comment}{\bfseries}{\itshape}
\spnewtheorem{fact}{Fact}{\bfseries}{\itshape}
  \baselineskip 8 mm

\newcommand{\comp}[1]{\overline{#1}}  

\begin{document}

\title{A Note on Threshold Dimension of Permutation Graphs}
\author {Diptendu Bhowmick
\thanks {Computer Science and Automation Department, Indian
Institute of Science, Bangalore- 560012 Email:
diptendubhowmick@gmail.com}}
\date{}
\institute{}
\maketitle

\begin{abstract}
A graph $G(V,E)$ is a threshold graph if there exist non-negative reals $w_v, \ v \in V$ and $t$ such that for every $U \subseteq V$, $\sum_{v \in U} w_v\leq t$ if and only if $U$ is a stable set. The {\it threshold dimension} of a graph $G(V,E)$, denoted as $t(G)$, is the smallest integer $k$ such that $E$ can be covered by $k$ threshold spanning subgraphs of $G$. A permutation graph is a graph that can be represented as the intersection graph of a family of line segments that connect two parallel lines in the Euclidean plane. In this paper we will show that if $G$ is a permutation graph then $t(G) \leq \alpha(G)$ (where $\alpha(G)$ is the cardinality of maximum independent set in $G$) and this bound is tight. As a corollary we will show that $t(G) \leq \frac{n}{2}$ where $n$ is the number of vertices in the permutation graph $G$. This bound is also tight.\\

\noindent {\bf Key words:} Threshold Dimension, Permutation Graph, Chromatic Number.
\end{abstract}

\section{Introduction}
Let $G$ be a simple, finite, undirected graph on $n$ vertices. The vertex set of $G$ is denoted as $V(G)$ and the edge set of $G$ is
denoted as $E(G)$. For any vertex $v \in V(G)$ let $N_G(v)=\{w \in V(G) \ | \ (v,w) \in E(G)\}$ be the set of neighbors of $v$. For each $S \subseteq V(G)$ let $G[S]$ denote the subgraph of $G$ induced by the vertices in $S$. Let $\alpha(G)$ and $\omega(G)$ denote the cardinality of maximum independent set and maximum clique in $G$ respectively. Also let $\chi(G)$ denote the chromatic number of $G$ i.e. $\chi(G)$ is the minimum number of colors needed for a proper vertex coloring of $G$.

Let $G'$ be a graph such that $V(G') = V(G)$. Then $G'$ is a {\it super graph} of $G$ if $E(G) \subseteq E(G')$. We define the {\it intersection} of two graphs as follows: if $G_1$ and $G_2$ are two graphs such that $V(G_1) = V(G_2)$, then the intersection of $G_1$ and $G_2$ denoted as $G = G_1 \cap G_2$ is a graph with $V(G) = V(G_1) = V(G_2)$ and $E(G) = E(G_1) \cap E(G_2)$.

A graph is a \textit{split graph} if its vertex set can be partitioned into a clique and an independent set. Split graphs were first studied by F\"{o}ldes and Hammer in \cite{aggregationInequalityChvatal,splitGraphFoldes}, and independently introduced by Tyshkevich and Chernyak \cite{canonicalPartitionTyshkevich}. For other characterizations and properties of split graphs one can refer to Golumbic \cite{algGraphTheoryPerfectGraphsGolumbic}.

 A graph $G(V,E)$ is a \textit{threshold graph} if there exist non-negative reals $w_v, \ v \in V$ and $t$ such that for every $U \subseteq V$, $\sum_{v \in U} w_v\leq t$ if and only if $U$ is a stable set. Let $2K_2,P_4,C_4$ denote a pair of independent edges, path on $4$ vertices and cycle on $4$ vertices respectively. Threshold graphs have a nice forbidden subgraph characterization as seen in the following theorem.

\begin{theorem}\textnormal{(Chv\'{a}tal and Hammer \cite{setPackingThresholdChvatalHammer})}
 A graph is a threshold graph if and only if it does not contain $2K_2,P_4$ or $C_4$ as induced subgraph.
\end{theorem}

 Chv\'{a}tal and Hammer \cite{setPackingThresholdChvatalHammer} introduced these graphs for their application in set-packing problems. In this paper we will use the following property of threshold graphs

\begin{fact}\textnormal{(see \cite{thresholdGraphsMahadevPeled} chapter 1)}\label{fact:tgDefn}
A graph $G$ is a threshold graph if and only if it is a split graph and for every pair of vertices $u,v$ in the independent set of $G$, either $N_G(u)\subseteq N_G(v)$ or $N_G(v)\subseteq N_G(u)$.
\end{fact}

\begin{definition}
A \emph{threshold cover} of a graph $G$ is a set of threshold graphs $\{T_1$ $,$ $T_2,$ $\cdots$ $,$ $T_k\}$ such that $V(G)=V(T_i)$ for $1 \leq i \leq k$ and $E(G)=E(T_1)\cup E(T_2)\cup \cdots \cup E(T_k)$. The \emph{threshold dimension} $t(G)$ is the least integer $k$ such that a threshold cover of cardinality $k$ exists for $G$. 
\end{definition}

\noindent{}Since complement of a threshold graph is also a threshold graph we have an equivalent definition of \emph{threshold dimension} as follows:
\begin{fact}
\emph{Threshold dimension} $t(G)$ of a graph $G$ is the smallest integer $k$ such that $\comp{G}$ can be represented as the intersection of $k$ threshold graphs. 
\end{fact}

Threshold graphs are the graphs having threshold dimension $1$. Chv\'{a}tal and Hammer \cite{aggregationInequalityChvatal} introduced the concept of threshold dimension. Ma \cite{phdThesisMa} has studied graphs having threshold dimension $2$. Recognition algorithms for graphs having threshold dimension $2$ were first proposed by Raschle and Simon \cite{recognitionRaschle} and improved by Sterbini and Raschle \cite{recognitionSterbini}. The problem of determining the threshold dimension of a graph has several applications like aggregation of linear inequalities in integer programming \cite{setPackingThresholdChvatalHammer,aggregationInequalityChvatal}, synchronization of competiting processes in a complex system such as a large computer \cite{henderson:88,Ordman85,Ordman86,Ordman89}, job scheduling \cite{Koop86}, Guttman scales in psychology \cite{thresholdDimCozzens,multiDimScalingCozzens} etc. 

Since every edge along with isolated vertices is a threshold subgraph, the threshold dimension is well-defined and is bounded by the number of edges of the graph. Chv\'{a}tal and Hammer \cite{aggregationInequalityChvatal} have shown the following upper bound on threshold dimension 
\begin{theorem}\textnormal{(Chv\'{a}tal and Hammer \cite{aggregationInequalityChvatal})}\label{chvatal}
 If $G$ is an undirected graph on $n$ vertices then $t(G) \leq n- \alpha(G)$. Moreover if $G$ is triangle-free then $t(G) = n- \alpha(G)$.
\end{theorem}

\noindent{}Chv\'{a}tal and Hammer have also shown that in general threshold dimension can be arbitrarily close to the number of vertices in the graph which is stated in the following theorem
\begin{corollary}\textnormal{(Chv\'{a}tal and Hammer \cite{setPackingThresholdChvatalHammer})}\label{close}
 For every $\epsilon>0$ there exists a graph $G$ with $n$ vertices such that $(1 - \epsilon)n < t(G)$.
\end{corollary}

\noindent{}Best known upper bound for threshold dimension of an $n$ vertex graph is shown by Erd\u{o}s et al. \cite{boundThresholdDimErdos}.
\begin{theorem}\textnormal{(Erd\u{o}s et al. \cite{boundThresholdDimErdos})}\label{erdos_bound}
 If $t(G)$ is the greatest threshold dimension of any graph $G$ on $n$ vertices then there exists some constant $A$ such that $n-A\sqrt{n}\log_2n<t(G)<n-\sqrt{n}+1$ where $n$ is large enough.
\end{theorem}

Chv\'{a}tal and Hammer \cite{aggregationInequalityChvatal} have shown that computation of threshold dimension for any graph is NP-hard. Yannakakis \cite{complexityPartialOrderDimnYannakakis} has shown that for any graph $G$ it is NP-complete to determine whether $t(G) \leq k$ where $k \geq 3$. Margot \cite{ComplexityMargot} has studied some complexity results about threshold graphs. For more references on threshold graphs and threshold dimension see the monograph of Mahadev Peled \cite{thresholdGraphsMahadevPeled}.

\subsection{Permutation Graphs:}
Let $\Pi$ be a permutation of the numbers $1,2, \ldots,n$. Then the graph $G[\Pi]=(V,E)$ is defined as follows: $V=\{1,2, \ldots, n \}$ and $(i,j) \in E \Leftrightarrow$ $(i-j)(\Pi^{-1}(i)-\Pi^{-1}(j))<0$, i.e. $i$ and $j$ occurs in the permutation in the reverse order. An undirected graph $G$ on $n$ vertices is called a permutation graph if there exists a permutation $\Pi$ of the numbers $1,2, \ldots,n$ such that $G \cong G[\Pi]$.

From the above definition it is easy to see that a permutation graph is the intersection graph of a family of line segments that connect two parallel lines in the Euclidean plane. We call such a family of line segments as the \textit{parallel line representation} of the permutation graph. Figure $1$ shows the permutation graph corresponding to the permutation $\Pi=\{4,7,5,1,2,6,3\}$ and its parallel line representation.
\begin{figure}[!ht]
\begin{minipage}{8cm}
\begin{tabular}{c}
\includegraphics[height=4cm]{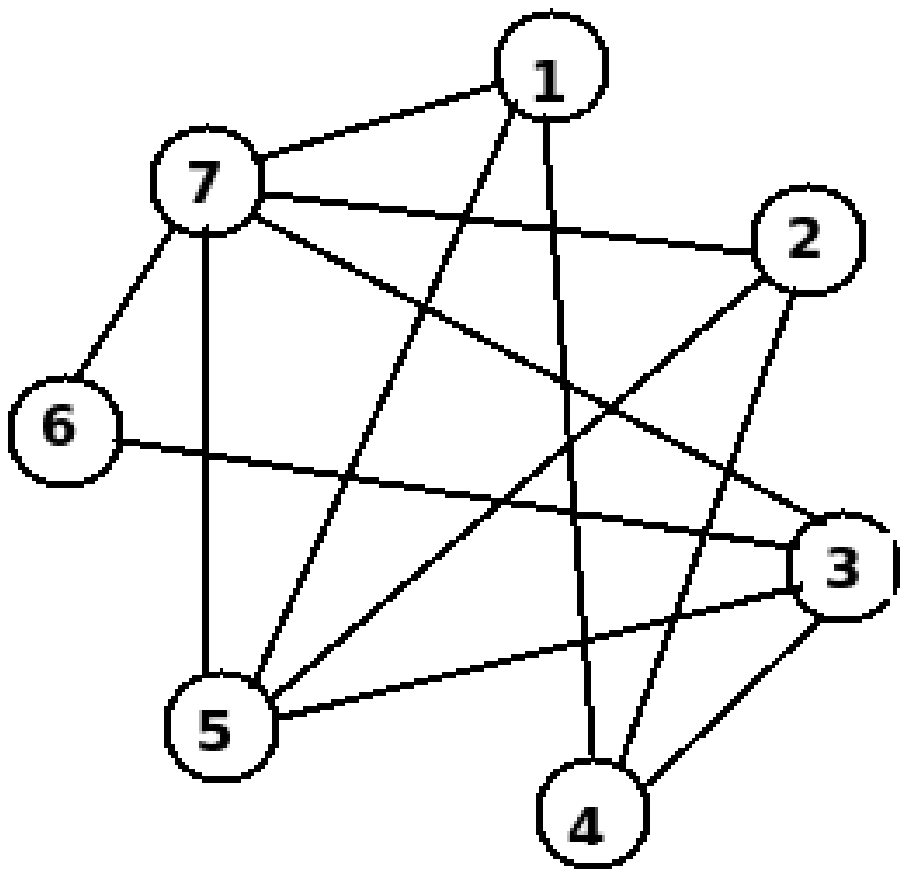}\\
(a)
\end{tabular}
\end{minipage}
\begin{minipage}{8cm}
\begin{tabular}{c}
\includegraphics[height=4cm]{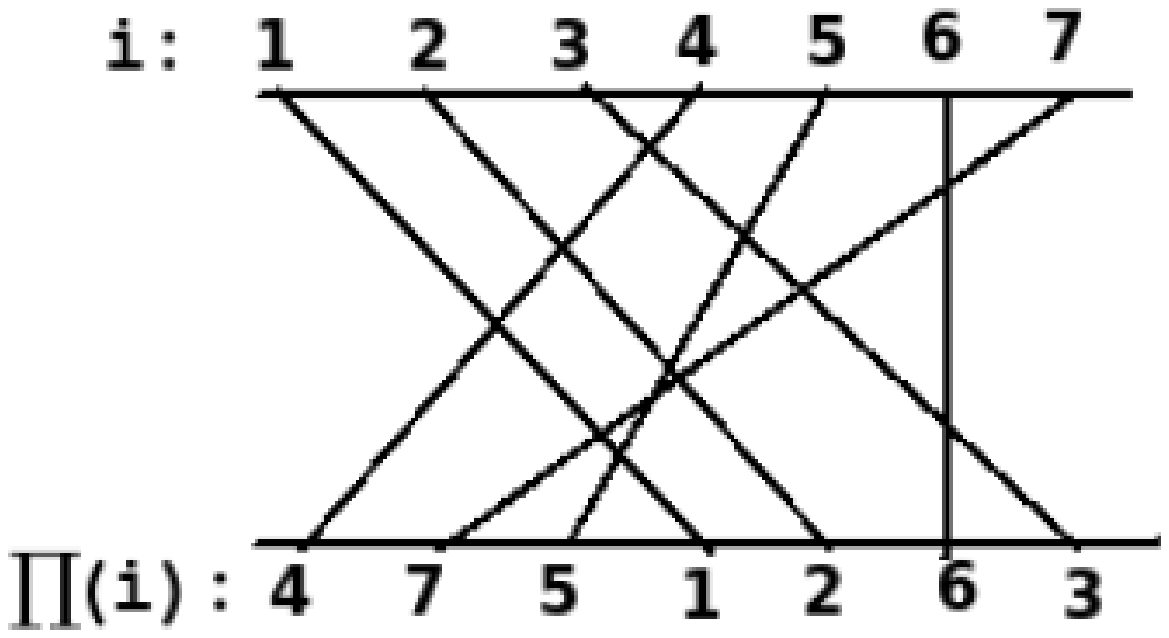}\\
(b)
\end{tabular}
\end{minipage}
\caption{An example. (a) A permutation graph and (b) its corresponding parallel line representation.}\label{fig:1}
\end{figure}
\begin{lemma}\label{Lem:ind_vertical}
 Let $G$ be a permutation graph corresponding to the permutation $\Pi$ and $I$ be an independent set of $G$. There exists a parallel line representation of $G$ which satisfies the following properties:
\begin{enumerate}
 \item Line segments corresponding to the vertices in $I$ are all vertical and distinct.
 \item Endpoints of the line segments corresponding to the vertices in $I$ come in the same order as their induced ordering in $\Pi$. 
\end{enumerate}
\end{lemma}
\begin{proof}
 Let $L$ be a parallel line representation of $G$. Since $I$ is an independent set the line segments corresponding to the vertices in $I$ are mutually non-intersecting in $L$. Therefore we can construct another parallel line representation $L'$ of $G$ by first placing the line segments corresponding to the vertices in $I$ so that they are all vertical and distinct and the ordering of their endpoints same as in $L$ (See how Figure \ref{fig:2} is obtained from Figure \ref{fig:1}). After this we can place the endpoints of the remaining line segments at distinct points on the lower and upper horizontal lines such that order of the endpoints on the lower as well as the upper horizontal lines are the same as in $L$. It is easy to see that $L'$ represents the graph $G$ and also satisfies the required properties.
\qed
\end{proof}

Figure (2) shows the parallel line representation of $G$ in which the line segments corresponding to the vertices of the independent set $\{1,2,3\}$ are all vertical.\\ 
\begin{figure}[!ht]
\centering
\includegraphics[height=4cm]{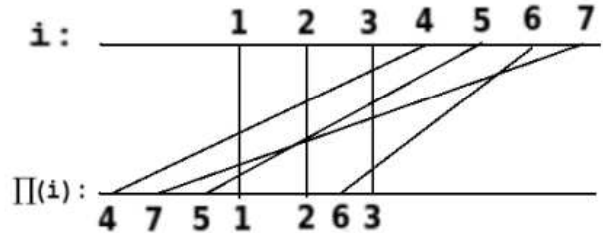}\\
\caption{Another representation of the example given in figure \ref{fig:1}}\label{fig:2}
\end{figure}

It is well known that permutation graphs are a subclass of perfect graphs. It is also a proper subclass of co-comparability graphs, comparability graphs and AT free graphs. 
\begin{fact}\label{fact:compl}
 An undirected graph $G$ is a permutation graph if and only if $G$ and $\overline{G}$ are comparability graphs and hence permutation graphs are closed under complementation.
\end{fact}

Permutation graphs are also a subclass of circle graphs (A circle graph is an intersection graph of chords in a circle). Permutation graph is a circle graph that admits an equator, i.e. one can draw an additional chord that intersects every other chord. See Golumbic \cite{algGraphTheoryPerfectGraphsGolumbic} for a brief introduction and references on permutation graphs.
\subsection{Our Results}
Let $G$ be a permutation graph on $n$ vertices. In this paper we will show that $t(G) \leq \alpha(G)$ and this bound is tight. As a corollary we will show that $t(G) \leq \frac{n}{2}$. This bound is also tight. Note that for general graphs both the bounds given above are not valid. Indeed as shown in Corollary \ref{close}, for general graphs threshold dimension can be arbitrarily close to $n$. The example given as proof of Corollary \ref{close} (See Golumbic \cite{algGraphTheoryPerfectGraphsGolumbic} chapter 10) illustrates that threshold dimension can be arbitrarily large compared to $\alpha(G)$. 
%
%
%
%
%
%
%

%
\section{Upper bound on threshold dimension of Permutation graphs}
In this section we will show that $t(G) \leq \alpha(G)$. Note that $\alpha(G)=\omega(\comp{G})=\chi(\comp{G})$, since $G$ is a permutation graph and hence a perfect graph. Therefore it suffices to show that $t(G) \leq \chi(\comp{G})$. Let $V(\comp{G})=\{1,2,\ldots,n\}$. According to Fact \ref{fact:compl}, $\comp{G}$ is also a permutation graph. Let $\Pi$ be the permutation corresponding to $\comp{G}$. So $(u,v) \in E(\comp{G})$, if and only if $(u-v)(\Pi^{-1}(u)-\Pi^{-1}(v))<0$. For each vertex $v \in V(\comp{G})$ let $\textbf{v}$ denote the line segment corresponding to the vertex $v$ in $L$, where $L$ is the parallel line representation of $\comp{G}$. Let $\chi(\comp{G})=k$ and $C_1,C_2,\ldots,C_{k}$ be the color classes corresponding to a proper vertex coloring of $\comp{G}$. 
\subsection{Index Set}
 Let $C_j=\{u_1,u_2,\cdots, u_{p(j)}\}$ where $p(j)=|C_j|$. Without loss of generality we can assume that $\Pi^{-1}(u_1)<\Pi^{-1}(u_2)<\cdots <\Pi^{-1}(u_p)$. The \textit{index set} $Ind_j(v)$ of a vertex $v$ $\in V(\comp{G}) \setminus C_j$ with respect to the set $C_j$ is the set of indices of vertices in $C_j$ to which $v$ is adjacent in $\comp{G}$ i.e. $Ind_j(v)=\{t:u_t \in C_j \ and \ (v,u_t) \in E(\comp{G}) \}$. For $v \in V(\comp{G}) \setminus C_j$, if $Ind_j(v) \neq \emptyset$ then \textit{minimum index} $l_j(v)$ of $v$ with respect to the set $C_j$ is defined to be $\min \ (Ind_j(v))$ and \textit{maximum index} $r_j(v)$ with respect  to the set $C_j$ is defined to be $\max \ (Ind_j(v))$.
\begin{lemma}\label{Lem:cont}
 For $v \in V(\comp{G}) \setminus C_j$, if $Ind_j(v) \neq \emptyset$ then $Ind_j(v)=\{t:l_j(v) \leq t \leq r_j(v)\}$.
\end{lemma}
\begin{proof}
 Let $l=l_j(v)$, $r=r_j(v)$. If $l=r$ then $|Ind_j(v)|=1$ and the Lemma is trivial. Therefore we assume that $l<r$. Since $C_j$ induces an independent set in $\comp{G}$, according to Lemma \ref{Lem:ind_vertical} we can construct a parallel line representation say $L$ of $\comp{G}$ such that the line segments corresponding to the vertices in $C_j$ are all vertical and distinct. Therefore ${\bf u_l}$ and ${\bf u_r}$ correspond to two distinct vertical lines in $L$. Since $l<r$, ${\bf u_l}$ lies to the left of ${\bf u_r}$ in $L$. Again since $\textbf{v}$ intersects both ${\bf u_l}$ and ${\bf u_r}$ and $l<r$, one endpoint of $\textbf{v}$ must lie to the left of ${\bf u_l}$ and other endpoint must lie to the right of ${\bf u_r}$. Thus $\textbf{v}$ intersects all the vertical line segments in between ${\bf u_l}$ and ${\bf u_r}$. Therefore $(v,u_t) \in E(\comp{G})$ for all $l \le t \le r$ and hence $Ind_j(v)=\{t:l_j(v) \leq t \leq r_j(v)\}$. 
\qed
\end{proof}
From Lemma \ref{Lem:cont}, it is clear that for each vertex $x=u_q \in C_j$ and $y \in V(\comp{G}) \setminus C_j$, if $(x,y) \notin E(\comp{G})$ and $Ind_j(y) \neq \emptyset$ then either $q<l_j(y)$ or $q>r_j(y)$. More specifically we have:

\begin{observation}\label{obs:1}
 Let $x=u_q \in C_j$ and $y \in V(\comp{G})\setminus C_j$. If $(x,y) \notin E(\comp{G})$ and $Ind_j(y) \neq \emptyset$ then either $q<l_j(y)$ or $q>r_j(y)$. Let $L_j$ be the parallel line representation of $G$ in which the line segments corresponding to the vertices in $C_j$ are all vertical and distinct.
\begin{enumerate}
 \item If $\Pi^{-1}(x) > \Pi^{-1}(y)$ then $q>r_j(y)$ and both endpoints of ${\bf y}$ lie to the left of (the vertical line segment) ${\bf x}$ in $L_j$.
 \item If $\Pi^{-1}(x) < \Pi^{-1}(y)$ then $q<l_j(y)$ and both endpoints of ${\bf y}$ lie to the right of (the vertical line segment) ${\bf x}$ in $L_j$.
\end{enumerate}
\end{observation}

\subsection{Threshold Graph Construction}
 We shall construct one threshold graph $T_j$ corresponding to each color class $C_j$ for $1 \le j \le k$ such that $\comp{G}=\bigcap_{j=1}^kT_j$. According to Lemma \ref{Lem:ind_vertical}, there exists a parallel line representation $L_j$ of $\comp{G}$ in which the line segments corresponding to the vertices in $C_j$ are all vertical and distinct. To construct $T_j$ we take the projection of each line segment on the lower horizontal line of $L_j$. Let $Proj_j(\textbf{v})$ denote the projection of $\textbf{v}$ on the lower horizontal line and $P$ denote the leftmost point of all the projections i.e. $P=\inf(\bigcup_{u\in V(\comp{G})}Proj_j(\textbf{u}))$. We map each vertex $v \in V(\comp{G})$ to an interval on the real line by the following mapping and define $T_j$ to be the intersection graph of the family of intervals $\{g_j(v):v \in V(\comp{G})\}$.

%
%
\begin{eqnarray*}
  	g_j(v) &=& [\inf Proj_j(\textbf{v}), \sup Proj_j(\textbf{v})] \ \ \ \ \ \ \ \ \ \ if \ v \in C_j\\
	&=& [P, \sup Proj_j(\textbf{v})] \ \ \ \ \ \ \ \ \ \ \ \ \ \ \ \ \ \ \ \ \ \ \ \ if \ v \in V(\comp{G}) \setminus C_j\\
\end{eqnarray*} 
\begin{comment}\label{Obs:distinct}
 For each $v \in C_j$, $g_j(v)$ corresponds to a distinct single point interval on the real line.
\end{comment}

\begin{lemma}
 $T_j$ is a threshold graph for $1 \le j \le k$.
\end{lemma}
\begin{proof}
From the construction it is easy to see that $T_j$ is an interval graph and the point $P$ is contained in all the intervals corresponding to the vertices in $V(\comp{G}) \setminus C_j$. Therefore $V(\comp{G}) \setminus C_j$ induces a clique in $T_j$. By Comment \ref{Obs:distinct}, the vertices in $C_j$ induce an independent set in $T_j$. Therefore $T_j$ is a split graph with $C_j$ as independent set and $V(\comp{G}) \setminus C_j$ as clique. Again since all the intervals corresponding to the vertices in $V(\comp{G}) \setminus C_j$ starts from the same point $P$ it is easy to see that for every two vertex $u,v \in C_j$ either $N_{T_j}(u) \subseteq N_{T_j}(v)$ or $N_{T_j}(v) \subseteq N_{T_j}(u)$. Therefore according to Fact \ref{fact:tgDefn}, $T_j$ is a threshold graph for $1 \le j \le k$.
\qed
\end{proof}

\begin{lemma}\label{Lem:supgraph}
For each threshold graph $T_j$, $E(\comp{G}) \subseteq E(T_j)$ for $1 \le j \le k$.
\end{lemma}
\begin{proof}
 Let $(x,y)\in E(\comp{G})$. Since $L_j$ is a parallel line representation of $\comp{G}$ there exists a point $Q\in \textbf{x} \cap \textbf{y}$ by definition. Clearly $Proj_j(Q) \in g_j(x) \cap g_j(y)$ and hence $g_j(x) \cap g_j(y) \ne \emptyset$. Therefore $E(\comp{G}) \subseteq E(T_j)$ for $1 \le j \le k$.
\qed
\end{proof}
\begin{lemma}\label{missing_edge}
If $(x, y) \notin E(\comp{G})$ then $(x, y) \notin E(T_j)$ for some $j$ where $1 \le j \le k$.
\end{lemma}
\begin{proof}
  Without loss of generality, assume that $x>y$. Let $x \in C_j$. Since $(x, y) \notin E(\comp{G})$ we have $\Pi^{-1}(x) > \Pi^{-1}(y)$. We consider the following cases:

\noindent{}\textbf{Case 1:} {\it When $y \in C_j$}. According to Comment \ref{Obs:distinct}, $g_j(x)$ and $g_j(y)$ correspond to two distinct points on the real line. Therefore $g_j(x) \cap g_j(y)=\emptyset$ and hence $(x,y) \notin E(T_j)$.\\

\noindent{}\textbf{Case 2:} {\it When $y \notin C_j$}. Let $x=u_q$ where $1 \leq q \leq |C_j|$. We consider the following cases:

\noindent{}\textbf{Subcase 2.1:} {\it When $Ind_j(y)=\emptyset$}. Clearly $(x, y) \notin E(\comp{G})$ and therefore ${\bf x}$ and ${\bf y}$ do not intersect in $L_j$. Since $(x,y) \notin E(\comp{G})$ and $\Pi^{-1}(x) > \Pi^{-1}(y)$, both endpoints of ${\bf y}$ lie to the left of ${\bf x}$ in $L_j$. Again since ${\bf x}$ is vertical in $L_j$, it corresponds to a point say $Q$ in $T_j$. Thus, $\sup Proj_j(\textbf{y}) <Q$. Therefore $g_j(x) \cap g_j(y)=\emptyset$ and hence $(x,y) \notin E(T_j)$.\\

 \noindent{}\textbf{Subcase 2.2:} {\it When $Ind_j(y) \neq \emptyset$}. Since $(x,y) \notin E(\comp{G})$ and $\Pi^{-1}(x) > \Pi^{-1}(y)$ according to Observation \ref{obs:1} part (1), we have $q>r_j(y)$ and both endpoints of $\textbf{y}$ lie to the left of (vertical line segment) ${\bf x=u_q}$ in $L_j$. Since ${\bf x}$ is vertical it corresponds to a point say $Q$ in $T_j$. Thus, $\sup Proj_j(\textbf{y}) <Q$. Therefore $g_j(x) \cap g_j(y)=\emptyset$ and hence $(x,y) \notin E(T_j)$.
\qed
\end{proof}
\noindent{}Combining Lemma \ref{Lem:supgraph} and \ref{missing_edge} we get $\bigcap_{j=1}^k T_j=\comp{G}$ and hence we have the following Theorem
\begin{theorem}\label{threshold_dim}
 If $G$ is a permutation graph then $t(G) \leq \chi(\comp{G})=\omega(\comp{G})= \alpha(G)$.
\end{theorem}
%


\noindent{} Combining Theorem \ref{chvatal} and Theorem \ref{threshold_dim} we get for a permutation graph, $t(G) \leq$ min($\alpha(G)$, $n-\alpha(G))$. As a corollary thus we have,

\begin{theorem}\label{thres_n/2}
 If $G$ is a permutation graph on $n$ vertices then $t(G) \leq \frac{n}{2}$.
\end{theorem}

\subsection{Tightness of Theorem \ref{threshold_dim} and \ref{thres_n/2}} 
\noindent{}\textbf{Example 1:} Let $G=(\frac{n}{2})K_2$ where $n$ is an even integer (i.e. a perfect matching on $n$ vertices). It is easy to see that $G$ is a permutation graph. Since $\alpha(G)=\frac{n}{2}$ we have $t(G) \leq \frac{n}{2}$ by Theorem \ref{threshold_dim}. It is also easy to see that $t(G) = \frac{n}{2}$ since it contains $\frac{n}{2}$ pairwise independent edges. So the upper bound for threshold dimension given in Theorem \ref{threshold_dim} and \ref{thres_n/2} is tight for $(\frac{n}{2})K_2$. \\

\noindent{}\textbf{Example 2:} Let $G=K_{\frac{n}{2},\frac{n}{2}}$ be a complete bipartite graph on $n$ vertices where $n$ is an even integer. It is easy to see that $G$ is a permutation graph. Thus we have $t(G) \leq \frac{n}{2}$ by Theorem \ref{threshold_dim}. But it was shown by Cozzens \cite{thresholdDimCozzens} that $t(G) = \frac{n}{2}$. So the upper bound for threshold dimension given in Theorem \ref{threshold_dim} and \ref{thres_n/2} is tight for $K_{\frac{n}{2},\frac{n}{2}}$. \\

\noindent{}\textbf{Example 3:} Let $G=\comp{(\frac{n}{2})K_2}$ where $n$ is an even integer (i.e. complement of a perfect matching on $n$ vertices). It is easy to see that $G$ is a permutation graph. Since $\alpha(G)=2$ we have $t(G) \leq 2$ by Theorem \ref{threshold_dim}. But it is easy to see that $G$ is not a threshold graph since it contains induced $C_4$. Therefore $t(G) \geq 2$. So the upper bound for threshold dimension given in Theorem \ref{threshold_dim} is tight for $\comp{(\frac{n}{2})K_2}$.

%
%
%
%
%
%

%
\end{document}